\documentclass[12pt]{article}
\usepackage{fullpage}
\usepackage{amssymb}
\usepackage{amsthm}
\usepackage[all,2cell]{xy}
\newtheorem{theorem}{Theorem}[section]
\newtheorem{corollary}[theorem]{Corollary}

\newtheorem{lemma}[theorem]{Lemma}

\newtheorem{proposition}[theorem]{Proposition}

\begin{document}
\title{Projective normality of quotient varieties modulo finite groups}
\author{S.S.Kannan, S.K.Pattanayak, Pranab Sardar  \\  
\\ Chennai Mathematical Institute, Plot H1, SIPCOT IT Park,\\ Padur 
Post Office, Siruseri, Tamilnadu - 603103, India.\\
kannan@cmi.ac.in, santosh@cmi.ac.in, pranab@cmi.ac.in} 

\maketitle
\date{}

\begin{abstract} In this note, we prove that for any finite dimensional vector 
space $V$ over an algebraically closed field $k$, and for any  finite subgroup 
$G$ of $GL(V)$ which is either solvable or is generated by pseudo reflections  
such that the $|G|$ is a unit in $k$, the projective variety $\mathbb P(V)/G$ 
is projectively normal with respect to the descent of 
$\mathcal O(1)^{\otimes |G|}$. 
\end{abstract}
\hspace*{4.5cm}Keywords: pseudo reflections, line bundle. 
\begin{center}
{\bf Introduction}
\end{center}

Let $G$ be a finite group. Let $V$ be a finite dimensional representation
of $G$ over a field $k$. In $1916$, E. Noether proved that if characterstic of $k$ not dividing $|G|$, then the $k$-algebra
of invariants $k[V]^G$ is finitely generated. In $1926$, she proved that the
same result holds in all characteristics. So, when $k$ is algebraically
closed, it is an interesting problem to study GIT- quotient varieties 
$V/G=Spec(k[V]^G)$ and $\mathbb P(V)/G$,(see [5] and [6]). Also, the line bundle
$\mathcal O(1)^{\otimes |G|}$ descends to the quotient $\mathbb P(V)/G$, where
$\mathcal O(1)$ denotes the ample generator of the Picard group of 
$\mathbb P(V)$. Let us denote it 
by $\mathcal L$. On the other hand, $V/G$ is normal. So, it is a natural
question to ask if $\mathbb P(V)/G$ is projectively normal with respect to the 
line bundle $\mathcal L$. In this note, we give an affirmative answer to this
question when:

\begin{enumerate}
\item $G$ is solvable and the characteristic of the base field is
not dividing $|G|$,

\item $G$ is a finite subgroup of $GL(V)$ generated by pseudo
reflections, where $V$ is a finite  dimensional vector space over a field  $k$
of characteristic not dividing $|G|$. 
\end{enumerate} 
Proof of the main result essentially uses an arithmetic result due to 
Erd\"{o}s-Ginzburg-Ziv (see [2]).
\section{Solvable case:}
In this section, we prove the following proposition which can be applied to
prove our main result when the group $G$ is solvable.
\begin{proposition}
Let $G$ be a finite solvable group, and let $V$ be a finite 
dimensional faithful representation of $G$ over a field $k$ of characteristic
not dividing $|G|$. Let $|G|=m$, 
$R:=\oplus_{d\geq 0}R_d$; $R_d:= (Sym^{dm}V)^G$. 
Then $R$ is generated as a $k$-algebra by $R_1$.
\end{proposition}

\begin{proof}
\underline{ Step $1$}: We first prove the statement when $G$ is 
cyclic of order $m$. Let $\xi$ be a $m$-th primitive root
of unity in an algebraic closure $\bar{k}$ of $k$. Let $F=k(\xi)$. 
Since $F$ is a free $k$ module, we have $V^G\otimes_k F=(V\otimes_k F)^G$. 
Hence, we may assume that $\xi \in k$. 

Let $G=<g>$.  Write $V=\oplus_{i=0}^{m-1}V_i$ 
where $V_i:=\{ v\in V: g.v=\xi^i.v\}$, $0\leq i\leq m-1$. Now let $f\in R_d$ be of the form $f=X_0\cdot X_1\cdots X_{m-1}$ with $X_i\in Sym^{a_i}V_i$ such that $ \sum_{i=0}^{m-1}a_i=dm$. Since $f$ is $G$-invariant we have 
\[ \sum_{i=0}^{m-1}i.a_i \equiv 0 \,\mbox{mod}\, m\]
If $d=1$, $f\in R_1$; so we may assume that $d\geq 2$. Now, consider the sequence of integers
\[ \underbrace{0,\ldots,0}_{a_0 \,\,\mbox{times}},\,\,\underbrace{1,\ldots,1}_{a_1\,\,\mbox{times}},\,\,\cdots, \,\, \underbrace{m-1,\ldots,m-1}_{a_{m-1} \,\,\mbox{times}}
\]
Since the sequence has $dm$ terms and $d\geq 2$, by a theorem of  Erd\"{o}s-Ginzburg-Ziv (see [2]), there is a subsequence with exactly $m$- terms whose terms add up to a
multiple of $m$. Thus there exist $f_1\in R_1$ and $f_2\in R_{d-1}$ such that $f=f_1.f_2$. Hence the proof follows by
induction on $deg(f)$. 

\underline{Step $2$}: Now we assume that $G$ is any finite solvable 
group of order $m$. We use induction on $m$ to prove the statement.
We may assume that $m$ is not a prime number.
Since $G$ is solvable it has a normal subgroup $H$ such that $G/H$
is a cyclic group of prime order. 

 Let $W:=(Sym^{|H|}V)^H$. Since $H$ is a normal subgroup of $G$,
both $\underbrace{ W\otimes \ldots \otimes W}_{d|G/H|\mbox{ copies}}$ and
$(Sym^{d|H|}V)^H$ have natural $G/H$-module structures. Let $G_1=G/H$.
Since $|H|<|G|$, by induction, the homomorphism 
$ \underbrace{ W\otimes \ldots \otimes W}_{d|G_1|\mbox{ copies}} \longrightarrow 
(Sym^{d|G|}V)^H$ is  surjective. $\hspace*{1cm}\cdots (1)$ 

{\bf Claim:} The natural map  $(Sym^{d|G_1|}W)^{G_1}\longrightarrow (Sym^{d|G|}V)^G$ is surjective.$\hspace{2cm}\ldots (2)$\\
The surjectivity of the natural map  $Sym^{d.|G_1|}W\longrightarrow (Sym^{d|G|}V)^H$ of $G_1$-modules follows from $(1)$ and  the following  commutative diagram
\[
\xymatrix{ \otimes^{d|G_1|}W \ar[d] \ar[r] & (Sym^{d|G|}V)^H  \\
       Sym^{d.|G_1|}W\ar[ru]    & }
\]
Hence applying Reynold's 
operator we have the claim.

 Now, consider the commutative diagram:
\[
\xymatrix{\otimes^d (Sym^{|G_1|}W)^{G_1} \ar[r] \ar[d] & (Sym^{d.|G_1|}W)^{G_1} \ar[d]\\
\otimes^d (Sym^{|G|}V)^{G} \ar[r] & (Sym^{d.|G|}V)^{G}}
\]
The first horizontal map  is surjective by step $(1)$ and the second 
vertical map is surjective by $(2)$. Thus the second horizontal map is
surjective. Thus the proposition follows.
\end{proof}

\section{Group generated by pseudo reflections:}

In this section, we will prove a combinatorial lemma which can be applied to 
prove our main result when the group $G$ is generated by pseudo reflections.

Let \underline{$a$}=$(a_1,a_2,\cdots a_r) \in \mathbb N^{r}$ and $N=(\prod_{i=1}^{r}a_{i})$. Consider the semigroup 

$M_{\underline{a}}=\{(m_1,m_2,\cdots m_r) \in \mathbb Z_{\geq 0}^{r}:\sum_{i=1}^{r}m_ia_i\equiv 0 \,\mbox{mod}\,N\}$ and the set 

$S_{\underline{a}}=\{(m_1,m_2,\cdots m_r) \in \mathbb Z_{\geq 0}^{r}:\sum_{i=1}^{r}m_ia_i= N\}$.

\begin{lemma}
 $M_{\underline{a}}$ is generated by $S_{\underline{a}}$ for \underline{$a$} $\in \mathbb N^{r}$.
\end{lemma}

\begin{proof}
Suppose $(m_1,m_2,\cdots m_r) \in \mathbb Z_{\geq 0}^{r}$ such that:
\[
\begin{array}{rcl}
\sum_{i=1}^{r}m_i.a_i=q.(\prod_{i=1}^{r}a_{i}), \mbox{with}\, \, q\geq 2.
\end{array}
\]
Let $a=\prod_{i=1}^{r}a_{i}$ and $n=q.a$. Then there exist a matrix \\
$A=$
$\pmatrix{x_{11} & x_{12} & \cdots & x_{1r}\cr x_{21} & x_{22} & \cdots & x_{2r}\cr \vdots & \vdots & \ddots & \vdots \cr   x_{n1} & x_{n2} & \cdots & x_{nr}}$   with $x_{ij}\in \mathbb Z_{\geq 0}$ for all $i$, $j$, such that $A\pmatrix{1 \cr 1 \cr \vdots \cr 1} = \pmatrix{1 \cr 1 \cr \vdots \cr 1}$, and $\pmatrix{1 \cr 1 \cr \vdots \cr 1} A = \pmatrix{ m_1 a_1 \cr m_2 a_2 \cr \vdots \cr m_r a_r}$.

Since $q\geq 2$, $n\geq 2.a_1$, the sequence $\{x_{11},x_{21},\cdots x_{n1}\}$
has atleast $2a_1$ number of terms. Hence applying the theorem of Erd\"{o}s-Ginzburg-Ziv (see [2]) and re arranging the entries of the first column if necessary we can assume that the $n$-terms of the sequence can be partitioned into $\frac{n}{a_1}$ number of subsequences $\{x_{11},x_{21},\cdots, x_{a_{1} 1}\}$, $\{x_{(a_1+1)1},x_{(a_1+2)1},\cdots, x_{2a_{1}1}\}$, $\cdots,$ $\{x_{(n-a_1+1)1},x_{(n-a_1+2)1},\cdots, x_{n1}\}$, each of length $a_1$ and sum of terms of each subsequence is a multiple of $a_1$. 

Again consider the sequence $\{\sum_{i=1}^{a_1}x_{i2}, \sum_{i=a_{1}+1}^{2a_1}x_{i2},\cdots \sum_{i=n-a_1+1}^{n}x_{i2}\}$. Using the same argument as above we can assume that this sequence can be partitioned into $\frac{n}{a_1a_2}$ number of subsequences each of length $a_2$ and sum of terms of each subsequence is a multiple of $a_2$.

Proceeding in this way, we can see that for each $j=1,2,\cdots r$, the sum of the first $a$-terms in the $j$th column of the matrix $A$ is a multiple of $a_j$.

 Let $b_j=\sum_{i=1}^{a}x_{ij}$. By construction of the $x_{ij}$'s , $b_{j}$ is a multiple of $a_{j}$ for every $j=1, 2, \cdots r$. So, for each $j=1, \cdots r$, there exists a non negative integer $b_{j}^{\prime}$ such that $b_{j}=a_{j}b_{j}^{\prime}$. Now, we see that the tuple $(b_{1}^{\prime},b_2^{\prime},\cdots b_r^{\prime})$ $ \in S_{\underline{a}}$, since 
$\sum_{j=1}^{r}b_j=a$. As $m_ja_{j}=\sum_{i=1}^{n}x_{ij}$, $b_{j}^{\prime}\leq m_{j}$ for every $j=1, 2, \cdots r$.
Thus, we have $(m_1,m_2,\cdots m_r)=(b_1^{\prime},b_2^{\prime},\cdots b_r^{\prime})+(m_1-b_1^{\prime},m_2-b_2^{\prime},\cdots m_r-b_r^{\prime})$, with $(b_1^{\prime},b_2^{\prime},\cdots b_r^{\prime})$ $ \in S_{\underline{a}}$. So, the lemma follows by induction on $q$.
\end{proof}

\begin{corollary}
Let $V$ be a finite dimensional vector space over a field $k$. Let $G$ be a
finite subgroup of GL$(V)$ which is generated by pseudo reflections. Further
assume that characteristic of $k$ is not dividing $|G|$, then $R=\oplus_{d\in\mathbb Z_{\geq 0}}R_d$ where
 $R_d:=(Sym^{d.|G| }(V^{*}))^G$, is generated by $R_1$.
\end{corollary}

\begin{proof}
By a theorem of Chevalley-Serre-Shephard-Todd (see [1],[4],[9],[10]),
$(Sym(V^{*}))^G$ is a polynomial ring $K[f_1,f_2,\cdots f_r]$ with each $f_i$ is a homogeneous polynomial of degree $d_i$ and $\prod_{i=1}^{r}d_i=|G|$. Thus, proof follows from lemma (2.1).
\end{proof}

We now prove our main result. 

Let $G$ be a finite group and $V$ be a finite 
dimensional, faithful representation of $G$ over an algebraically closed 
field of characteristic not dividing $|G|$. Let $\mathcal O(1)$ denote the 
ample generator of the Picard group of $\mathbb P(V)$. Let $\mathcal L$ 
denote the descent of the line bundle $\mathcal O(1)^{\otimes |G|}$ to the quotient $\mathbb P(V)/G$. Further, 
assume that either $G$ is solvable or is  generated by pseudo reflections 
in $GL(V)$. Then, we have:

\begin{theorem}
$\mathbb P(V)/G$ is projectively normal with respect to $\mathcal L$.
\end{theorem} 

\begin{proof}
Proof follows from Proposition (1.1), Corollary (2.2), and exercise $5.14(d)$, Chapter II of (3).
\end{proof}

\end{document}